\documentclass[a4paper, leqno, 12pt]{article}
\usepackage{amsmath,amsthm}
\usepackage{amssymb,latexsym}
\usepackage{enumerate}
\usepackage[T1]{fontenc}
\usepackage{array}
\usepackage[all]{xypic}

\usepackage[left=3cm,right=3cm,top=3.0cm,bottom=3.0cm]{geometry}

\newtheorem{thm}{Theorem}[section]

\newtheorem{lm}[thm]{Lemma}

\theoremstyle{definition}
\newtheorem{df}[thm]{Definition}

\numberwithin{equation}{section}

 1

\def\QQ{\mathbb{Q}}

\def\NN{{\mathbb{N}}}

\def\CL{{\cal L}}

\def\CP{{\cal P}}

\def\CS{{\cal S}}

\def\epv {{$\mbox{}$\hfill ${\Box}$\vspace*{1.5ex} }}

\def\Gen{\textnormal{Gen}}

\def\ra{\rightarrow}

\def\ov{\overline}
\def\mod{\mbox{{\rm mod}}}
\def\lmod{\mbox{{\rm -mod}}}

\def\End{\mbox{{\rm End}}}

\def\ncong{{\cong\hspace{-3mm} / \hspace{1mm} }}

\def\ov#1{\overline{#1}}

\begin{document}


\baselineskip=17pt


\title{On wild algebras and super-decomposable pure-injective modules}

\author{Grzegorz Pastuszak${}^{*}$}

\date{}

\maketitle

\renewcommand{\thefootnote}{}

\footnote{${}^{*}$Faculty of Mathematics and Computer Science, Nicolaus Copernicus University, Chopina 12/18, 87-100 Toru\'n, Poland, past@mat.uni.torun.pl.}

\begin{abstract}
Assume that $k$ is an algebraically closed field and $A$ is a finite-dimensional wild $k$-algebra. Recently, L. Gregory and M. Prest proved that in this case the width of the lattice of all pointed $A$-modules is undefined. Hence the result of M. Ziegler implies that there exists a super-decomposable pure-injective $A$-module, if the base field $k$ is countable. Here we give a straightforward proof of the fact that there exists a special family of pointed $A$-modules, called an independent pair of dense chains of pointed $A$-modules. This also yields the existence of a super-decomposable pure-injective $A$-module.
\end{abstract}

\section{Introduction}


The remarkable tame and wild dichotomy of Yu. Drozd \cite{Dr} states that the class of finite-dimensional algebras over algebraically closed fields divides into two disjoint classes: \textit{tame algebras} and \textit{wild algebras}. The class of wild algebras properly contains the class of \textit{strictly wild algebras}. We refer the reader to \cite{SiSk3} for definitions of these classes. Furthermore, A. Skowro\'nski introduced in \cite{SkBC} a concept of the \textit{growth} of a tame algebra. This yields a stratification of the class of tame algebras into \textit{domestic}, \textit{linear} and \textit{polynomial growth} algebras. Tame algebras which are not of polynomial growth are called \textit{non-polynomial growth} algebras. 

Understanding various aspects of representation types is still one of the central topics of the representation theory of finite-dimensional algebras over algebraically closed fields. A good example supporting this fact is provided by the tame \textit{self-injective} algebras. Indeed, the representation theory of these algebras is well developed for the polynomial growth (see \cite{Sk0,Sk4}), but much less is known for the non-polynomial growth. Recently, K. Erdmann and A. Skowro\'nski introduced in \cite{ErSk5} (and study in a huge ongoing project, see the introduction of \cite{ErSk6} for more details) a prominent class of \textit{weighted surface algebras}. These algebras are some special representation-infinite tame symmetric algebras (and hence self-injective). Since most of them is of non-polynomial growth, they play a significant role in understanding the representation theory of all tame self-injective algebras.

The representation type is studied by using various concepts and methods. On the level of finite-dimensional modules we have, in particular, results on the shape of connected components of the Auslander-Reiten quiver or the component quiver, see for example \cite{Sk96}, \cite{Sk98} or \cite{JawSk}. On the level of infinite-dimensional modules, a fundamental characterization of the representation type is given in \cite{CB} in terms of \textit{generic modules}. The paper \cite{GKM} studies representation type in terms of matrix problems and matrix reduction algorithms. It is conjectured by M. Prest (see \cite{Pr}, \cite{Pr22}) that a finite-dimensional algebra $A$ over an algebraically closed field is of domestic representation type if and only if the \textit{Krull-Gabriel dimension} $\textnormal{KG}(A)$ of $A$ is finite (see \cite{Ge1,Ge2} for definitions and \cite{P4} for a list of results supporting this conjecture). The second conjecture due to Prest states that an algebra $A$ is of domestic representation type if and only if there is no \textit{super-decomposable pure-injective} $A$-module (see for example \cite{Pr2}). Such a module is some special infinite-dimensional module. This paper is related to the second conjecture of Prest.

Assume that $R$ is a ring with a unit. By an \textit{$R$-module} we mean a left $R$-module. An $R$-module $M$ is \textit{super-decomposable} if and only if $M\neq 0$ and $M$ has no indecomposable direct summands. We refer to \cite{Kiel}, \cite{HZ} and \cite[Chapter 7]{JeLe} for the concept of \textit{pure-injectivity}.

The problem of the existence of super-decomposable pure-injective $R$-modules is stated in \cite{Zi}. In this paper M. Ziegler proves a fundamental criterion for such modules to exist, asserting that if the ring $R$ is countable, then $R$ possesses a super-decomposable pure-injective module if and only if the width of the lattice of all \textit{pp-formulas} is undefined, see \cite{Pr} or \cite{Pr22} for the definitions. The later statement can be formulated in terms of the lattice of all \textit{pointed finitely-presented $R$-modules} (these lattices are isomorphic, see \cite{PPT} and \cite{KaPa2} for more details). 

The case when $R$ is a finite-dimensional algebra over a field $k$ is studied in many papers. We refer to the introduction of \cite{KaPa3} for an up-to-date list of results in this direction (except for the most recent one, stating that representation-infinite domestic standard self-injective algebras over algebraically closed fields do not have super-decomposable pure-injective modules, see \cite[Theorem 8.3]{P4}). In most of these papers it is assumed that the base field $k$ is countable. This yields $R$ is countable and hence one can apply the criterion of Ziegler. All the known results support the conjecture of Prest concerning super-decomposable pure-injective modules.

The first result on the existence of super-decomposable pure-injective modules for finite-dimensional algebras is proved by M. Prest in \cite[Theorem 13.7]{Pr}. It states that these modules do exist over strictly wild algebras. For a very long time it was not known whether this holds for wild algebras. Recently, L. Gregory and M. Prest prove in \cite{GP} that this is the case. Indeed, they show in \cite[Theorem 2.1]{GP} that any representation embedding functor induces an embedding of lattices of pp-formulas. This implies that the width of the lattice of all pp-formulas over a wild algebra $A$ is undefined, and so there exists a super-decomposable pure-injective $A$-module, if the base field is countable (see Corollary 2.4 of \cite{GP}). We stress that the paper \cite{GP} contains many more interesting results.

Recall that an \textit{independent pair of dense chains of pointed modules} is some special family of pointed modules that allows to formulate a handy sufficient condition for the existence of a super-decomposable pure-injective module (see Section 2 for the details). This notion is introduced in \cite{PPT} and generalized in \cite{KaPa2}. It is successfully used in the series of papers \cite{KaPa,KaPa2,Pa,KaPa3} which is devoted to the problem of the existence of super-decomposable pure-injective modules for \textit{strongly simply connected algebras}, see \cite{Sk2,NoSk}. 

This paper is devoted to show that if $A$ is a wild $k$-algebra over an algebraically closed field $k$, then there exists an independent pair of dense chains of pointed $A$-modules. This result is proved in Theorem 4.2. The proof of Theorem 4.2 is rather straightforward. It is based on some facts from \cite{KaPa} on the existence of independent pairs of dense chains of pointed modules for string algebras of non-polynomial growth and the very definition of a wild algebra. 

Let us clarify the overlap of \cite{GP} and the present paper. It is proved in Corollary 2.4 of \cite{GP} that the width of the lattice of all pp-formulas over a wild algebra is undefined. We prove in Theorem 4.2 the existence of independent pairs of dense chains of pointed modules for wild algebras. Thus Theorem 4.2 implies \cite[Corollary 2.4]{GP} (see Theorem 2.4), but the converse is not known (however unexpected, see Theorem 2.3). In this sense, Theorem 4.2 is stronger than \cite[Corollary 2.4]{GP}. However, we observe that, after some additional work, the assertion of Theorem 4.2 could be derived from \cite[Theorem 2.1]{GP}. Nevertheless, the paper \cite{GP} of Gregory and Prest takes a different perspective than the point of view presented here. Therefore we see our results as another, simple and independent of \cite{GP}, proof of the fact that wild algebras possess super-decomposable pure-injective modules.

The paper contains four sections. In Section 2 we collect basic information on pointed modules, recall the definition of an independent pair of dense chains of pointed modules and formulate the Ziegler criterion. We also introduce some special independent pairs of dense chains of pointed modules which we call \textit{strong}. Section 3 is devoted to show the existence of a strong independent pair of dense chains of pointed modules (with some additional properties) over any string algebra of non-polynomial growth, see Theorem 3.4. This result is only implicitly contained in \cite{KaPa}, so we decided to include here its detailed proof. This makes the present paper more convenient to the reader. In Section 4 we present our main results. Indeed, Theorem 4.1 shows that representation embedding functors, under some mild assumptions, preserve strong independent pairs of dense chains of pointed modules. Then Theorem 4.2, stating the existence of an independent pair of dense chains of pointed modules and a super-decomposable pure-injective module for any wild algebra (if the base field is countable), is a direct consequence of Theorem 4.1 and Theorem 3.4. 

Throughout, $k$ is a fixed algebraically closed field. By an \textit{algebra} we mean a finite-dimensional associative basic $k$-algebra with a unit. If $A$ is an algebra, then by an \textit{$A$-module} we mean a left $A$-module. We denote by $A\lmod$ the category of all finitely-generated (hence finite-dimensional) left $A$-modules.




\section{The lattice of pointed modules and a sufficient existence condition}

In this section we recall some basic facts on pointed modules and related concepts. In particular, we present a sufficient condition for the existence of a super-decomposable pure-injective module in terms of independent pairs of dense chains of pointed modules.

Assume that $R$ is a ring with a unit. We denote by $R\lmod$ the category of all finitely-presented left $R$-modules. Assume that $\Theta\in R\lmod$. A \textit{$\Theta$-pointed $R$-module} is a pair $(M,\chi_{M})$ where $M$ is a finitely-presented $R$-module and $\chi_{M}:\Theta\ra M$ is an $R$-module homomorphism.

Assume that $(M,\chi_{M})$ and $(N,\chi_{N})$ are $\Theta$-pointed $R$-modules. By a \textit{$\Theta$-pointed $R$-homomorphism} from $(M,\chi_{M})$ to $(N,\chi_{N})$ we mean an $R$-homomorphism $f:M\ra N$ such that $f\chi_{M}=\chi_{N}$. If $f:M\ra N$ is a $\Theta$-pointed $R$-homomorphism from $(M,\chi_{M})$ to $(N,\chi_{N})$, we write $f:(M,\chi_{M})\ra(N,\chi_{N})$. If $f:M\ra N$ is an isomorphism, we call $f:(M,\chi_{M})\ra(N,\chi_{N})$ a \textit{$\Theta$-pointed isomorphism} and the corresponding $\Theta$-pointed modules $(M,\chi_{M})$ and $(N,\chi_{N})$ \textit{$\Theta$-isomorphic}.

Assume that $t\in\NN$, $t\geq 1$, $\Theta=R^{t}$ and $(M,\chi_{M})$ is a $\Theta$-pointed $R$-module. Assume that $e_{1},\hdots,e_{t}$ form the $R$-base of the module $\Theta$. The homomorphism $\chi_{M}$ is uniquely determined by the elements $\chi(e_{1}),\hdots,\chi(e_{t})\in M$. This yields that any $\Theta$-pointed $R$-module can be identified with a tuple $(M,m_{1},\hdots,m_{t})$ where $M$ is an $R$-module and $m_{1},\hdots,m_{t}\in M$. Moreover, a $\Theta$-pointed $R$-homomorphism from $(M,m_{1},\hdots,m_{t})$ to $(N,n_{1},\hdots,n_{t})$ can be identified with an $R$-homomorphism $f:M\ra N$ such that $f(m_{i})=n_{i}$, for $i=1,\hdots,t$. In case $\Theta=R$, we simply speak about pointed modules and pointed homomorphisms.

Let $P_{R}^{\Theta}$ be the set of all $\Theta$-isomorphism classes of $\Theta$-pointed $R$-modules. Let $\equiv$ be a binary relation on $P_{R}^{\Theta}$ defined by $(M,\chi_{M})\equiv(N,\chi_{N})$ if and only if there exist pointed homomorphisms $f:(M,\chi_{M})\ra(N,\chi_{N})$ and $g:(N,\chi_{N})\ra(M,\chi_{M})$. Then $\equiv$ is an equivalence relation and the quotient set ${\cal{P}}_{R}^{\Theta}=P_{R}^{\Theta}\slash\equiv$ is a poset with respect to the relation $\leq$ defined by $\ov{(M,\chi_{M})}\leq\ov{(N,\chi_{N})}$ if and only if there exists a pointed homomorphism $f:(N,\chi_{N})\ra(M,\chi_{M})$. We denote by $\ov{(S,\chi_{S})}$ the $\equiv$-class of a $\Theta$-pointed $R$-module $(S,\chi_{S})$.

The poset ${\cal{P}}_{R}^{\Theta}$ is a modular lattice with respect to the operations $\oplus$ and $*$ defined below, see \cite{Pr} for details. 

Assume that $(M,\chi_{M})$, $(N,\chi_{N})$ are $\Theta$-pointed $R$-modules. A $\Theta$-pointed $R$-module $(M\oplus N,\chi_{M\oplus N})$ where $\chi_{M\oplus N}(l)=(\chi_{M}(l),\chi_{N}(l))$ for any $l\in\Theta$ is the \textit{pointed direct sum} of $(M,\chi_{M})$ and $(N,\chi_{N})$. We set $(M,\chi_{M})\oplus(N,\chi_{N})=(M\oplus N,\chi_{M\oplus N})$.

Assume that $M*N$ is the pushout of $\chi_{M}$ and $\chi_{N}$, that is, $$M*N=M\oplus N/\{(\chi_{M}(l),-\chi_{N}(l));l\in\Theta\}.$$ Moreover, let $\epsilon_{M}:M\ra M*N$, $\epsilon_{N}:N\ra M*N$ be the $R$-module homomorphisms given by $\epsilon_{M}(m)=\overline{(m,0)}$, $\epsilon_{N}(n)=\overline{(0,n)}$ for any $m\in M$, $n\in N$. A $\Theta$-pointed $R$-module $(M*N,\chi_{M*N})$ where $\chi_{M*N}=\epsilon_{M}\chi_{M}=\epsilon_{N}\chi_{N}$ is the \textit{pointed pushout} of $(M,\chi_{M})$ and $(N,\chi_{N})$. We set $(M,\chi_{M})*(N,\chi_{N})=(M*N,\chi_{M*N})$.

It is easy to see that $$\sup\{\ov{(M,\chi_{M})},\ov{(N,\chi_{N})}\}=\ov{(M\oplus N,\chi_{M\oplus N})},$$$$\inf\{\ov{(M,\chi_{M})},\ov{(N,\chi_{N})}\}=\ov{(M*N,\chi_{M*N})}.$$ Recall that if $\Theta=R^{t}$, then the lattice ${\cal{P}}_{R}^{\Theta}$ is equivalent to the lattice of all \textit{pp-formulas} with $t$ free variables ($t\geq 1$), see \cite{Pr,Pr22}. 

We recall definitions of wide lattices of pointed modules and independent pairs of dense chains of pointed modules. Moreover, we present in Theorem 2.3 the relation between these notions.

We say that a lattice $\CL\subseteq{\cal{P}}_{R}^{\Theta}$ of $\Theta$-pointed $R$-modules is \textit{wide} if and only if for any $\ov{(M_{p},\chi_{M_{p}})}<\ov{(M_{q},\chi_{M_{q}})}\in\CL$ there are incomparable elements $\ov{(M,\chi_{M})},\ov{(N,\chi_{N})}$ of $\CL$ such that $$\ov{(M_{p},\chi_{M_{p}})}<\ov{(M,\chi_{M})},\ov{(N,\chi_{N})}<\ov{(M_{q},\chi_{M_{q}})},$$ $$\ov{(M_{p},\chi_{M_{p}})}\leq\ov{(M*N,\chi_{M*N})}<\ov{(M\oplus N,\chi_{M\oplus N})}\leq\ov{(M_{q},\chi_{M_{q}})}.$$ In case the lattice ${\cal{P}}_{R}^{\Theta}$ contains a wide sublattice $\CL$, we say that the \textit{width of ${\cal{P}}_{R}^{\Theta}$ is undefined}. The above definition is a special case of a general definition of a wide lattice, see \cite{Pr,Pr22} or Section 3 of \cite{KaPa2}.

Assume that $C$ is a set. A family $\{(M_{q},\chi_{M_{q}});q\in C\}$ of $\Theta$-pointed $R$-modules is denoted by $(M_{q},\chi_{M_{q}})_{q\in C}$. Let $\QQ$ be the set of rational numbers viewed as a poset with respect to the natural ordering $\leq$. Recall that a poset $P$ is a \textit{$\QQ$-chain} if and only if it is a dense chain without end points. It is well known that any $\QQ$-chain is isomorphic as a poset with the set $\QQ$. 

Assume that $R$ is a ring with a unit and $\Theta$ is a finitely-presented $R$-module. The following definitions were introduced in \cite{PPT} and generalized in \cite{KaPa2}. 

\begin{df} Assume that $C$ is a $\QQ$-chain. A \textit{dense chain of $\Theta$-pointed $R$-modules} is a family $(M_{q},\chi_{M_{q}})_{q\in C}$ of $\Theta$-pointed $R$-modules such that:

{\rm (a)} the endomorphism ring $\End_{R}(M_{q})$ is local and $\chi_{M_{q}}\neq 0$ for any $q\in C$,

{\rm (b)} there exist $\Theta$-pointed homomorphisms $\mu_{q,q'}:(M_{q},\chi_{M_{q}})\ra(M_{q'},\chi_{M_{q'}})$ for any $q<q'\in C$,  

{\rm (c)} the pointed modules $(M_{q},\chi_{M_{q}})$ and $(M_{q'},\chi_{M_{q'}})$ are not $\Theta$-isomorphic for any $q\neq q'\in C$. \epv
\end{df}

\begin{df} An \textit{independent pair of dense chains of $\Theta$-pointed $R$-modules} is a pair $((M_{q},\chi_{M_{q}})_{q\in C_{1}},(N_{t},\chi_{N_{t}})_{t\in C_{2}})$ of dense chains of $\Theta$-pointed $R$-modules such that:

{\rm (a)} the endomorphism ring $\End_{R}(M_{q}*N_{t})$ is local for any $q\in C_{1},t\in C_{2}$ where $(M_{q}*N_{t},\chi_{M_{q}*N_{t}})=(M_{q},\chi_{M_{q}})*(N_{t},\chi_{N_{t}})$,

{\rm (b)} the pointed module $(M_{q},\chi_{M_{q}})*(N_{t},\chi_{N_{t}})$ is not $\Theta$-isomorphic to $(M_{q'},\chi_{M_{q'}})*(N_{t},\chi_{N_{t}})$ nor to $(M_{q},\chi_{M_{q}})*(N_{t'},\chi_{N_{t'}})$ for any $q\neq q'\in C_{1},t\neq t'\in C_{2}$. \epv
\end{df}

Independent pairs of dense chains of pointed modules generate wide lattices of pointed modules in the following way.

\begin{thm} Assume that the pair $((M_{q},\chi_{M_{q}})_{q\in C_{1}},(N_{t},\chi_{N_{t}})_{t\in C_{2}})$ is an independent pair of dense chains of $\Theta$-pointed $R$-modules. Then the lattice $$\Gen(\ov{(M_{q},\chi_{M_{q}})}_{q\in C_{1}}\cup\ov{(N_{t},\chi_{N_{t}})}_{t\in C_{2}}),$$which is the smallest sublattice of ${\cal{P}}_{R}^{\Theta}$ containing sets $\ov{(M_{q},\chi_{M_{q}})}_{q\in C_{1}}$ and $\ov{(N_{t}\chi_{N_{t}})}_{t\in C_{2}}$, is a wide lattice. Therefore the width of the lattice ${\cal{P}}_{R}^{\Theta}$ is undefined.
\end{thm}

{\bf Proof.} The assertion is a direct consequence of \cite[Theorem 3.4]{KaPa2}. \epv

It is not known whether the existence of a wide sublattice of ${\cal{P}}_{R}^{\Theta}$ (or, equivalently, a wide sublattice of the lattice of all pp-formulas over $R$) implies the existence of an independent pair of dense chains of $\Theta$-pointed $R$-modules.

The assertion $(1)$ of the following theorem is the Ziegler's criterion, see \cite{Zi}, and $(2)$ is a handy version of this criterion. Observe that $(2)$ follows directly from $(1)$ and Theorem  2.3.

\begin{thm} Assume that $R$ is a countable ring with a unit and $\Theta$ is a finitely presented $R$-module. 
\begin{enumerate}[\rm(1)]
	\item If the lattice ${\cal{P}}_{R}^{\Theta}$ of $\Theta$-pointed $R$-modules has width undefined, then there exists a super-decomposable pure-injective $R$-module.
	\item If there exists an independent pair of dense chains of $\Theta$-pointed $R$-modules, then there exists a super-decomposable pure-injective $R$-module. \epv
\end{enumerate}
\end{thm}

We apply the above theorem only when $R$ is a finite-dimensional $k$-algebra. Note that in this case any finitely-presented $R$-module $M$ is finite-dimensional and $\End_{R}(M)$ is local if and only if $M$ is indecomposable. Moreover, it is easy to see that $R$ is countable if and only if the field $k$ is countable. 

Assume that $(M,\chi_{M})$ and $(N,\chi_{N})$ are $\Theta$-pointed $R$-modules. Note that if we have $(M,\chi_{M})\cong(N,\chi_{N})$, then $M\cong N$, but the converse does not hold in general. This justifies the following special version of Definitions 2.1 and 2.2.

\begin{df} Assume that $R$ is a ring with a unit. 
\begin{enumerate}[\rm(1)]
	\item A dense chain $(M_{q},\chi_{M_{q}})_{q\in C}$ of $\Theta$-pointed $R$-modules is \textit{strong} if and only if $M_{q}$ and $M_{q'}$ are not isomorphic (as $R$-modules), for any $q<q'\in C$.
	\item An independent pair $((M_{q},\chi_{M_{q}})_{q\in C_{1}},(N_{t},\chi_{N_{t}})_{t\in C_{2}})$ of dense chains of $\Theta$-pointed $R$-modules is \textit{strong} if and only if dense chains $(M_{q},\chi_{M_{q}})_{q\in C_{1}}$ and $(N_{t},\chi_{N_{t}})_{t\in C_{2}}$ are strong and the module $M_{q}*N_{t}$ is not isomorphic (as an $R$-module) to $M_{q'}*N_{t}$ nor to $M_{q}*N_{t'}$, for any $q\neq q'\in C_{1},t\neq t'\in C_{2}$. \epv
\end{enumerate}
\end{df}

\section{Independent pairs of dense chains for string algebras}

In this section we recall Theorem 5.7 from \cite{KaPa} (see Theorem 3.3) and derive its important special case (see Theorem 3.4). This theorem (which is some refinement of \cite[Theorem 4.1]{Pu1}) states that if $A$ is a string algebra of non-polynomial growth, then there exists and independent pair of dense chains of pointed $A$-modules. We stress that this pair is strong, see Definition 2.5. This fact does not play a significant role in \cite[Theorem 5.7]{KaPa} itself, but is a crucial ingredient in proofs of our main results. 
 
The content of this section is faithfully based on Sections 4 and 5 of \cite{KaPa}. For convenience, some concepts related with string algebras are recalled.

Let $Q=(Q_{0},Q_{1})$ be a finite quiver with the set $Q_{0}$ of vertices and the set $Q_{1}$ of arrows. Given an arrow $\alpha\in Q_{1}$ with the starting point $s(\alpha)$ and the terminal point $t(\alpha)$ we denote by $\alpha^{-1}$ its \textit{formal inverse}. We set $s(\alpha^{-1})=t(\alpha)$, $t(\alpha^{-1})=s(\alpha)$ and $(\alpha^{-1})^{-1}=\alpha$. The set of all formal inverses of the arrows from $Q_{1}$ is denoted by $Q_{1}^{-1}$. The elements of $Q_{1}$ are called \textit{direct arrows} whereas of $Q_{1}^{-1}$-\textit{inverse arrows}.

By a \textit{walk} from $x$ to $y$ of length $n\geq 1$ in $Q$ we mean a sequence $c_{1}\hdots c_{n}$ in $Q_{1}\cup Q_{1}^{-1}$ such that $s(c_{n})=x\in Q_{0}$, $t(c_{1})=y\in Q_{0}$, $s(c_{i})=t(c_{i+1})$ and $c_{i}^{-1}\neq c_{i+1}$, for all $1\leq i<n$. We agree that $(c_{1}\hdots c_{n})^{-1}=c_{n}^{-1}\hdots c_{1}^{-1}$. A walk $c_{1}\hdots c_{n}$ is a \textit{path} provided $c_{i}\in Q_{1}$, for $1\leq i\leq n$. Furthermore, to each vertex $x\in Q_{0}$ we associate the \textit{stationary path} $e_{x}$ of length $0$, with $s(e_{x})=t(e_{x})=x$.

Given a finite quiver $Q=(Q_{0},Q_{1})$ we denote by $kQ$ the \textit{path algebra} of the quiver $Q$. The $k$-basis of $kQ$ is the set of all paths in $Q$ and the multiplication in $kQ$ is induced by the concatenation of paths. For example, if $\alpha,\beta\in Q_{1}$ and $s(\alpha)=t(\beta)$, then $\alpha\beta$ is the path $\stackrel{\beta}{\longrightarrow}\stackrel{\alpha}{\longrightarrow}$.

A two-sided ideal $I$ in $kQ$ is called \textit{admissible} if $\langle Q_{1}\rangle^{n}\subseteq I\subseteq\langle Q_{1}\rangle^{2}$, for some $n\in\NN$, $n\geq 2$. If $I$ is an admissible ideal in $kQ$, then the pair $(Q,I)$ is called the \textit{bound quiver} and the associated quotient algebra $kQ\slash I$ the \textit{bound quiver algebra}. The fundamental result of P. Gabriel \cite{Ga0,Ga01} states that any finite-dimensional basic associative $k$-algebra over algebraically closed field $k$ is isomorphic to some bound quiver $k$-algebra, see also Chapter II of \cite{AsSiSk}.

A bound quiver $(Q,I)$ and the corresponding bound quiver $k$-algebra $kQ\slash I$ is \textit{special biserial} \cite{SkWa} if an only if the following conditions are satisfied:

\begin{itemize}
	\item any vertex of $Q$ is the starting point of at most two arrows and the terminal point of at most two arrows,
	\item given an arrow $\beta$ there is at most one arrow $\alpha$ with $s(\beta)=t(\alpha)$ and $\beta\alpha\notin I$ and at most one arrow $\gamma$ with $s(\gamma)=t(\beta)$ and $\gamma\beta\notin I$.
\end{itemize}

A \textit{string algebra} is a special biserial algebra $kQ/I$ such that $I$ is generated by paths. By a \textit{string} in the string algebra $kQ/I$ we mean a walk $c_{1}\hdots c_{n}$ in $Q$ such that neither $c_{i}\hdots c_{i+t}$ nor $c_{i+t}^{-1}\hdots c_{i}^{-1}$ belongs to $I$, for $1\leq i<i+t\leq n$. Moreover, by a \textit{band} in $kQ/I$ we mean a string $S=c_{1}\hdots c_{n}$ such that:

\begin{itemize}
	\item all powers of $S$ are defined, i.e. $t(c_{1})=s(c_{n})$ and $S^{m}$ is a string, for all $m\in\NN$, 
	\item $c_{1}$ is a direct arrow and $c_{n}$ is an inverse arrow.
\end{itemize}

Assume that $S=c_{1}\hdots c_{n}$ is a string in the string algebra $kQ/I$. We denote by $M(S)$ the associated \textit{string module}. This is some $(n+1)$-dimensional module with $k$-linear basis $\{z_{1},\hdots,z_{n}\}$, called the \textit{canonical basis} of $M(S)$ and denoted as a tuple $(z_{1},\hdots,z_{n})$. Recall that any string module is indecomposable and $M(S_{1})\cong M(S_{2})$ if and only if $S_{1}=S_{2}$ or $S_{1}=S_{2}^{-1}$. We refer to \cite{BuRi} and \cite{SkWa} and for more details on string algebras and string modules.

Assume that $A=kQ\slash I$ is a fixed string algebra. Given an arrow $a\in Q_{1}$, we define $\CS(a)$ to be the set of the strings over $A$ that start with $a$. We recall from \cite{BuRi} that there exists some linear ordering on $\CS(a)$. We denote this ordering by $\leq$ and set $S<T$ if and only if $S\leq T$ and $S\neq T$.

A pair $(U,V)$ of two different bands over $A$ starting with the same direct arrow and ending with the same inverse arrow is {\em $\QQ$-generating} provided $U<V$ and $U$ is not a prolongation of $V$ or vice versa. This means that  $U\neq VX$ and $V\neq UY$, for any strings $X$ and $Y$. Moreover, assume that $\Sigma(U,V)$ is the set of all finite words over the alphabet $\{U,V\}$, including the empty word $\phi$.

\begin{thm}\textnormal{\cite[Theorem 5.3]{KaPa}} Assume that $(U,V)$ is a $\QQ$-generating pair of bands over $A$ and $S,T\in\Sigma(U,V)$. Then the set $\CL_{S}^{T}(U,V)=\{SXTU;X\in\Sigma(U,V)\}$ is a dense chain without end points. \epv
\end{thm}

Assume that $S=s_{1}\hdots s_{n}$ is a string over the string algebra $A$, $M(S)$ is the associated string module and $z_{1}^{S}\in M(S)$ is the first element of the canonical $k$-basis of $M(S)$. We call the pointed $A$-module $(M(S),z_{1}^{S})$ the \textit{canonical pointed string module} associated with $S$.

Assume that $T$, $S$ are strings over $A$ such that $T=t_{1}\hdots t_{k}$, $S=s_{1}\hdots s_{m}$ and $TS$ is also a string. We denote by $z^{(T,S)}$ the element $z_{k+1}^{TS}$ of the canonical basis $(z_{1}^{TS},\hdots ,z_{k+1}^{TS},\hdots ,z_{k+m+1}^{TS})$ of $M(TS)$.

The following fact is proved in \cite[3.1,3.2]{Pu2}.

\begin{lm} Asume that $A=kQ\slash I$ is a string algebra.
\begin{enumerate}[\rm(1)]
	\item Assume that $a\in Q_{1}$, $S,T\in\CS(a)$ and $S<T$. There exists a pointed $A$-homomorphism $f_{(T,S)}:(M(T),z_{1}^{T})\ra(M(S),z_{1}^{S})$ of the canonical pointed string modules $(M(T),z_{1}^{T})$ and $(M(S),z_{1}^{S})$. 
	\item Assume that $T,S$ are strings over $A$ such that $T^{-1}S$ is also a string. The pointed module $(M(T^{-1}S),z^{(T^{-1},S)})$ is the pointed pushout of the pointed modules $(M(S),z_{1}^{S})$ and $(M(T),z_{1}^{T})$. \epv
\end{enumerate}
\end{lm}

The following theorem is proved in \cite{KaPa}. Here we stress the fact that $\QQ$-generating pairs of bands over string algebras induce independent pairs of dense chains of pointed modules that are strong.

\begin{thm}\textnormal{\cite[Theorem 5.7]{KaPa}} Assume that $(U,V)$ and $(U^{-1},V^{-1})$ are $\QQ$-generating pairs of bands over the string algebra $A$. Let $S,T\in\Sigma(U,V)$ and $S',T'\in\Sigma(U^{-1},V^{-1})$. Then the pair $$((M(X),z_{1}^{X})_{X\in\CL_{S}^{T}(U,V)},(M(Y),z_{1}^{Y})_{Y\in\CL_{S'}^{T'}(U^{-1},V^{-1})})$$ is a strong independent pair of dense chains of pointed modules in $A\lmod$. 
\end{thm}

{\bf Proof.} Set $C_{1}=\CL_{S}^{T}(U,V)$, $C_{2}=\CL_{S'}^{T'}(U^{-1},V^{-1})$. We prove that $(M(X),z_{1}^{X})_{X\in C_{1}}$ and $(M(Y),z_{1}^{Y})_{Y\in C_{2}}$ are strong dense chains of pointed $A$-modules. First observe that Theorem 3.1 yields the sets $C_{1},C_{2}$ are dense chains without end points. 

The modules $M(X)$ and $M(Y)$ are indecomposable, for any strings $X\in C_{1}$ and $Y\in C_{2}$, since they are string modules over $A$.

The existence of pointed homomorphisms in $(M(X),z_{1}^{X})_{X\in C_{1}}$, $(M(Y),z_{1}^{Y})_{Y\in C_{2}}$ follows from Lemma 3.2 (1).

The modules $M(X_{1})$ and $M(X_{2})$ are not isomorphic, for any $X_{1},X_{2}\in C_{1}$ such that $X_{1}\neq X_{2}$. Indeed, $X_{1}\neq X_{2}^{-1}$ since $X_{2}^{-1}$ starts with a different direct arrow than $X_{1}$ and $X_{1}\neq X_{2}$ by the assumption. Similarly, $M(Y_{1})\ncong M(Y_{2})$, for any $Y_{1},Y_{2}\in C_{2}$ such that $Y_{1}\neq Y_{2}$.

Consequently, $(M(X),z_{1}^{X})_{X\in C_{1}}$ and $(M(Y),z_{1}^{Y})_{Y\in C_{2}}$ are strong dense chains of pointed modules in $A$-$\mod$. Now we prove that they form a strong independent pair.

Assume that $X,X_{1},X_{2}\in C_{1}$, $Y,Y_{1},Y_{2}\in C_{2}$. The pointed pushout of $(M(X),z_{1}^{X})$ and $(M(Y),z_{1}^{Y})$ is indecomposable since it is isomorphic with $(M(X^{-1}Y),z^{(X^{-1},Y)})$, by Lemma 3.2 (2).

If $X_{1}\neq X_{2}$, then $M(X_{1}^{-1}Y)\ncong M(X_{2}^{-1}Y)$, because $X_{1}^{-1}Y\neq X_{2}^{-1}Y$ (since $X_{1}\neq X_{2}$) and $X_{1}^{-1}Y\neq (X_{2}^{-1}Y)^{-1}=Y^{-1}X_{2}$ (since $X_{1}^{-1}$ starts with a different direct arrow than $Y^{-1}$). Similarly, if $Y_{1}\neq Y_{2}$, then $M(X^{-1}Y_{1})\ncong M(X^{-1}Y_{2})$.

The above arguments show that $((M(X),z_{1}^{X})_{X\in C_{1}},(M(Y),z_{1}^{Y})_{Y\in C_{2}})$ is a strong independent pair of dense chains of pointed $A$-modules. \epv

We recall that if the base field $k$ is algebraically closed and $A$ is a string algebra, then $A$ is tame of non-polynomial growth if and only if there exist $\QQ$-generating pairs of bands $(U,V)$, $(U^{-1},V^{-1})$, see \cite{ri,ri2,Sch,WaWa}. 

In Section 4 we apply the following refinement of Theorem 3.3. 

\begin{thm} There exists a string algebra $\Lambda$ and a strong independent pair $\CP=((M_{q},\chi_{M_{q}})_{q\in C_{1}},(N_{t},\chi_{N_{t}})_{t\in C_{2}})$ of dense chains of $\Theta$-pointed $\Lambda$-modules such that $\Theta$ is indecomposable and $\chi_{M_{q}},\chi_{N_{t}}$ are monomorphisms, for any $q\in C_{1}$, $t\in C_{2}$. 
\end{thm}

{\bf Proof.} We give a concrete example of a string algebra satisfying the thesis of the theorem. This algebra plays an important role in \cite{KaPa} (see also \cite{KaPa2} and \cite{Pa}). Assume that \begin{center}$\xymatrix{\\ Q=}$$\xymatrix{x_{1} \ar@/^0.3pc/[d]^{\beta} \ar@/_0.3pc/[d]_{\alpha}\\x_{2} \ar@/^0.3pc/[d]^{\gamma}\ar@/_0.3pc/[d]_{\delta}\\x_{3}}$\end{center} and $\Lambda=kQ\slash I$, where $I=\langle\delta\alpha,\gamma\beta\rangle$. It is easy to see that $\Lambda$ is a string algebra. Moreover, direct calculations show that if $(U,V)=(\gamma\alpha\beta^{-1}\delta^{-1},\gamma\delta^{-1})$, then $(U,V)$ and $(U^{-1},V^{-1})$ are $\QQ$-generating pairs of bands over $\Lambda$. 

We set $C_{1}=\CL^{\phi}_{\phi}(U,V)$ and $C_{2}=\CL^{\phi}_{\phi}(U^{-1},V^{-1})$. Then Theorem 3.3 yields that the pair $$\CP_{1}=((M(X),z_{1}^{X})_{X\in C_{1}},(M(Y),z_{1}^{Y})_{Y\in C_{2}})$$ is a strong independent pair of dense chains of pointed modules in $\Lambda\lmod$. 

Let $P(x_{3})=\Lambda e_{x_{3}}$ be the simple projective $\Lambda$-module associated to the vertex $x_{3}$. We recall from \cite{KaPa} (see particularly Corollary 7.1 (a)) that $z_{1}^{X}\in e_{x_{3}}M(X)$ and $z_{1}^{Y}\in e_{x_{3}}M(Y)$, for any $X\in C_{1}$, $Y\in C_{2}$. Observe that if $(M,m)$ is a pointed $\Lambda$-module such that $m\in e_{x_{3}}M$, then there is a homomorphism $\chi_{M}:P(x_{3})\ra M$ such that $\chi_{M}(e_{x_{3}})=m$, so $(M,m)$ can be identified with $(M,\chi_{M})$. These arguments imply that $\CP_{1}$ induces a strong independent pair $$\CP_{2}=((M(X),\chi_{M(X)})_{X\in C_{1}},(M(Y),\chi_{M(Y)})_{Y\in C_{2}})$$ of dense chains of $\Theta$-pointed $\Lambda$-modules where $\Theta=P(x_{3})$ (see Lemma 3.10 of \cite{KaPa2} for the details). Observe that $\Theta$ is a simple module, so homomorphisms $\chi_{M(X)},\chi_{M(Y)}$ are monomorphisms, for any $X\in C_{1}$, $Y\in C_{2}$. Therefore, the thesis holds for the pair $\CP=\CP_{2}$. \epv

\section{The main results}

This section is devoted to prove our main results. In Theorem 4.1 we show that representation embeddings preserve strong independent pairs of dense chains of pointed modules which additionally satisfy conditions from the thesis of Theorem 3.4. It follows directly from Theorem 4.1 and Theorem 3.4 that any wild algebra $A$ possesses an independent pair of dense chains of pointed modules. Hence there exists a super-decomposable pure-injective $A$-module, if the base field is countable. These facts are stated in Theorem 4.2.

Throughout the section, $A,B$ are $k$-algebras. All functors considered are covariant functors. Assume that $F:B\lmod\ra A\lmod$ is a functor and $(M,\chi_{M})$ is a $\Theta$-pointed $B$-module, for some $\Theta\in B\lmod$. The $F(\Theta)$-pointed $A$-module $(F(M),F(\chi_{M}))$ is denoted by $F(M,\chi_{M})$. If $(M_{q},\chi_{M_{q}})_{q\in C}$ is a family of $\Theta$-pointed $B$-modules, then $F(M_{q},\chi_{M_{q}})_{q\in C}$ denotes the family $(F(M_{q}),F(\chi_{M_{q}}))_{q\in C}$ of $F(\Theta)$-pointed $A$-modules.

Recall that a functor $F:B\lmod\ra A\lmod$ is a \textit{representation embedding} if and only if $F$ is exact, respects the isomorphism classes (that is, $F(X)\cong F(Y)$ implies $X\cong Y$, for any $B$-modules $X,Y$) and carries indecomposable modules to indecomposable ones. An algebra $A$ is of \textit{wild representation type} (or \textit{wild}) if and only if there exists a representation embedding functor $F:C\lmod\ra A\lmod$, for any $k$-algebra $C$.

\begin{thm} Assume that $A,B$ are $k$-algebras and $F:B\lmod\ra A\lmod$ is a representation embedding. Assume that $((M_{q},\chi_{M_{q}})_{q\in C_{1}},(N_{t},\chi_{N_{t}})_{t\in C_{2}})$ is a strong independent pair of dense chains of $\Theta$-pointed $B$-modules such that the module $\Theta$ is indecomposable and homomorphisms $\chi_{M_{q}},\chi_{N_{t}}$ are monomorphisms, for any $q\in C_{1}$, $t\in C_{2}$. Then $$(F(M_{q},\chi_{M_{q}})_{q\in C_{1}},F(N_{t},\chi_{N_{t}})_{t\in C_{2}})$$ is an independent pair of dense chains of $F(\Theta)$-pointed $A$-modules.
\end{thm} 

{\bf Proof.} We show that $F(M_{q},\chi_{M_{q}})_{q\in C_{1}}$ is a dense chain of $F(\Theta)$-pointed $A$-modules (similar arguments show that $F(N_{t},\chi_{N_{t}})_{t\in C_{2}}$ is a dense chain of $F(\Theta)$-pointed $A$-modules as well). Indeed, assume that $q\in C_{1}$. The module $F(M_{q})$ is indecomposable, because $M_{q}$ is indecomposable. Since $F(\chi_{M_{q}}):F(\Theta)\ra F(M_{q})$ is a monomorphism such that $F(\Theta)$ indecomposable, we get $F(\chi_{M_{q}})\neq 0$. 

Assume that $q<q'$ and $\mu_{q,q'}:(M_{q},\chi_{M_{q}})\ra(M_{q'},\chi_{M_{q'}})$ is a $\Theta$-pointed homomorphism. Since $\mu_{q,q'}\chi_{M_{q}}=\chi_{M_{q'}}$, we get $F(\mu_{q,q'})F(\chi_{M_{q}})=F(\chi_{M_{q'}})$, so $F(\mu_{q,q'})$ is a $F(\Theta)$-pointed homomorphism from $F(M_{q},\chi_{M_{q}})$ to $F(M_{q'},\chi_{M_{q'}})$. 

Furthermore, $F(\Theta)$-pointed modules $F(M_{q},m_{q})$ and $F(M_{q'},m_{q'})$ are not isomorphic, for any $q\neq q'$. Indeed, we have $F(M_{q})\ncong F(M_{q'})$, because $M_{q}\ncong M_{q'}$. This shows that $F(M_{q},\chi_{M_{q}})_{q\in C_{1}}$ is a dense chain of $F(\Theta)$-pointed $A$-modules.

We show that that dense chains $F(M_{q},\chi_{M_{q}})_{q\in C_{1}}$ and $F(N_{t},\chi_{N_{t}})_{t\in C_{2}}$ form an independent pair. Indeed, the module $F(M_{q}*N_{t})$ is indecomposable, because $M_{q}*N_{t}$ is indecomposable, for any $q\in C_{1}, t\in C_{2}$. 

Observe that the functor $F:\mod(B)\ra\mod(A)$ preserves finite colimits (and finite limits), since it is exact, see \cite{McL}. This implies that $$F((M_{q},\chi_{M_{q}})*(N_{t},\chi_{N_{t}}))\cong F(M_{q},\chi_{M_{q}})*F(N_{t},\chi_{N_{t}}),$$ for any $q\in C_{1},t\in C_{2}$. In particular, we get $F(M_{q}*N_{1})\cong F(M_{q})*F({N_{t}})$. Assume that $F(M_{q},\chi_{M_{q}})*F(N_{t},\chi_{N_{t}})$ is isomorphic with $F(M_{q},\chi_{M_{q}})*F(N_{t'},\chi_{N_{t'}})$, for some $q\in C_{1}$ and $t\neq t'\in C_{2}$. Then we get $$F(M_{q}*N_{t})\cong F(M_{q})*F(N_{t})\cong F(M_{q})*F(N_{t'})\cong F(M_{q}*N_{t'}),$$ which yields $M_{q}*N_{t}\cong M_{q}*N_{t'}$. Since this is not the case, we get that $F(M_{q},\chi_{M_{q}})*F(N_{t},\chi_{N_{t}})$ is not isomorphic with $F(M_{q},\chi_{M_{q}})*F(N_{t'},\chi_{N_{t'}})$. Similar arguments show that it is not isomorphic with  $F(M_{q'},\chi_{M_{q'}})*F(N_{t},\chi_{N_{t}})$ as well, for any $q'\neq q\in C_{1}$. This shows the assertion. \epv

\begin{thm} Assume that $A$ is a wild $k$-algebra over an algebraically closed field $k$. There exists an independent pair of dense chains of $\Xi$-pointed modules, for some $A$-module $\Xi$. Therefore there exists a super-decomposable pure-injective $A$-module, if the base field $k$ is countable.
\end{thm}

{\bf Proof.} Since $A$ is a wild algebra, there exists a representation embedding functor $F:\Lambda\lmod\ra A\lmod$ where $\Lambda$ is the string algebra considered in Theorem 3.4. Thus the assersion follows from Theorem 3.4, Theorem 4.1 and Ziegler's criterion (Theorem 2.4 (2)). \epv

\end{document}